\documentclass[11pt,a4paper]{article}

\usepackage{amsmath}
\usepackage{amsthm}
\usepackage{amssymb}

\numberwithin{equation}{section} \allowdisplaybreaks

\newtheorem{proposition}{Proposition}[section]
\newtheorem{corollary}{Corollary}[section]

\theoremstyle{definition}
\newtheorem{definition}{Definition}[section]

\newtheorem{remark}{Remark}[section]

\begin{document}
\thispagestyle{empty}
\begin{center}
\end{center}
\bigskip

\centerline{{\sc On the modular classes of Poisson-Nijenhuis manifolds}}

\bigskip

\centerline{Yvette Kosmann-Schwarzbach and Franco Magri}

\bigskip

\noindent{\it Abstract}
We prove a property of the
Poisson-Nijenhuis manifolds which yields new proofs of the bihamiltonian
properties of the hierarchy of modular
vector fields defined by Damianou and Fernandes.
\bigskip

\bigskip

\section*{Introduction}

In \cite{DF}, Damianou and
Fernandes  
defined the modular vector field and the modular
class of a Poisson-Nijenhuis manifold, 
and they proved that the hierarchy generated by the modular
vector field coincides with the canonical hierarchy of bihamiltonian vector
fields already defined in \cite{M}. 
A theorem of Beltr\'an and Monterde \cite{BM} states that, in a
PN-manifold, the derived bracket (see e.g. \cite{K})
of the interior products by $N$ and
$P$ acting on forms is the interior product by the hamiltonian
vector field with hamiltonian $- \frac{1}{2} {\rm{Tr}} N$. 
In this Letter, 
we give an elementary 
proof of a particular case of this theorem,
a simple consequence of which, stated in 
Corollary \ref{corollary}, 
enables us to give new proofs of the hamiltonian 
properties of the hierarchy of modular
vector fields of PN-manifolds.
These can be extended to the case of arbitrary PN-algebroids in a
straightforward manner.

\bigskip

\section{Poisson-Nijenhuis structures}
There are many ways of expressing the compatibility of a pair $(P,N)$,
where $N$ is a
Nijenhuis tensor and $P$ is  
a Poisson bivector on a manifold $M$ satisfying the condition that 
$NP$ be skew symmetric, in order to ensure that $NP, N^2P, 
\ldots, N^kP, \ldots$ be a sequence of pairwise-compatible Poisson
brackets. 
Let $d_N= [i_N,d]$ be the differential on forms 
associated with the deformed bracket of vector fields, 
$[ \, , \, ]_N$, and let $[ \, , \, ]_P$ be the graded bracket of forms
defined by $P$. When no confusion is possible, 
we denote by 
$N$ both the Nijenhuis tensor and its transpose, and by $P$ both the
Poisson bivector and the map from 1-forms to vectors it defines, with
the convention $P\alpha =i_{\alpha}P$.
Let $H^P_f= Pdf$ be the hamiltonian vector field with hamiltonian $f \in
C^{\infty}(M)$ in the Poisson structure $P$.
The derived bracket $[[i_N,d],i_P]= [d_N,i_P]$ is denoted by $[i_N,i_P]_d$.

\newpage

\begin{proposition} 
The following conditions on $N$ and $P$ are equivalent:
\begin{itemize}
\item
(i) $NP=P N$ and
(ii) $C(P,N) = 0$, where, for all $\alpha$, $\beta \in \Gamma (T^*M)$,
$$C(P,N)(\alpha, \beta) = [\alpha, \beta]_{NP} - ([N\alpha, \beta]_P +
[\alpha, N\beta]_P - N [\alpha, \beta]_P) \ .
$$ 
\item
$d_N$ is a derivation of bracket
$[ \, , \, ]_P$. 
\item 
$d_P = [P, \cdot ]$ is a derivation of the deformed bracket
$[ \, , \, ]_N$. 
\item 
Let $\{ \, , \, \}_{NP}$ be 
the Poisson bracket of functions with respect to $NP$.

(i) $NP=P N$ and
(ii) $d \{ f , g \}_{NP} = L_{H^P_f}d_N g -L_{H^P_g}d_N f - d_N(H^P_f( g))$,
for all $f$, $g \in C^{\infty}(M)$.
\end{itemize}
\end{proposition}
This last condition follows from $C(P,N)(df, dg) = 0$, for all functions 
$f$, $g \in C^{\infty}(M)$, 
using the relation $[\alpha, df]_P = - i_{H^P_f} d\alpha$.

\begin{definition} When any one of the above conditions
is satisfied,
$N$ and $P$ are called {\it compatible}. 
The pair $(P,N)$ is a {\it Poisson-Nijenhuis structure} (or
PN-structure) if
$N$ and $P$ are compatible. A manifold with a Poisson-Nijenhuis
structure is called a {\it Poisson-Nijenhuis manifold} (or PN-manifold).
\end{definition}

The compatibility of $P$ and $N$ can also be stated in terms of the
morphism properties of maps $P$, $N^k P$, $N^k$ and $(^tN)^k$, 
$k \geq 1$, relating the
various Lie algebroid structures on $TM$ and $T^*M$. 

\begin{proposition}\label{prop1}
Let $P$ be a Poisson bivector and $N$ a Nijenhuis tensor on $M$ 
such that $PN=NP$.
Then, for all $\alpha \in \Gamma (T^*M)$,
\begin{equation}\label{trace}
\frac{1}{2}  {\rm{Tr}}  (C(P,N)\alpha) = \frac{1}{2} <Pd {\rm{Tr}}  \, N, \alpha > +
[i_N,i_P]_d \alpha \ ,
\end{equation}
where $[ \, ,\, ]_d$ denotes the derived bracket.
\end{proposition} 
\begin{proof}
We shall use the expression of the components of $C(P,N)$ in local
coordinates \cite{PN},
$$
C^{kj}_m =  
P^{lj}_{} \partial_{l}  N^{k}_{m} + P^{kl}_{} \partial_{l} N^{j}_{m} 
- N^{l}_{m} \partial_{l}  P^{kj}_{} + N^{j}_{l} \partial_{m}  P^{kl}_{}
- P^{lj}_{} \partial_{m}  N^{k}_{l} \ ,
$$ 
whence 
$$
C^{kj}_k =  
P^{lj} \partial_{l}  N^{k}_{k} + P^{kl}_{} \partial_{l} N^{j}_{k} 
- N^{l}_{k} \partial_{l}  P^{kj}_{} + N^{j}_{l} \partial_{k}  P^{kl}_{}
- P^{lj}_{} \partial_{k}  N^{k}_{l} \ .
$$ 
From the assumption $NP=PN$, i.e., $P^{lj}N^k_l+ P^{lk}N^j_l =0$,
we obtain
$$
N^{k}_{l} \partial_{m}  P^{lj}_{} + P^{lj}_{} \partial_{m} N^{k}_{l} 
+ N^{j}_{l} \partial_{m}  P^{lk}_{} 
+ P^{lk}_{} \partial_{m}  N^{j}_{l} = 0 \ ,
$$
whence 
$$
N^{k}_{l} \partial_{k}  P^{lj}_{} + P^{lj}_{} \partial_{k} N^{k}_{l} 
+ N^{j}_{l} \partial_{k}  P^{lk}_{} 
+ P^{lk}_{} \partial_{k}  N^{j}_{l} = 0 \ .
$$
This identity implies that
$$
\frac{1}{2} C^{kj}_k = \frac{1}{2} P^{lj}\partial_l N^k_k + P^{lk}
\partial_k N^j_l \ .
$$
Thus, for any 1-form $\alpha$,
$$
\frac{1}{2}  {\rm{Tr}}  (C(P,N)\alpha) = 
\frac{1}{2} P^{lj}\partial_l N^k_k \alpha_j + P^{lk}
\partial_k N^j_l \alpha_j
$$
$$ =  - \frac{1}{2} < Pd {\rm{Tr}} N, \alpha >
 + i_P d i_N \alpha - i_{NP} d\alpha   \ .
$$
Since $i_{NP} = i_{PN} = i_{P}i_{N}$,
$$
(i_P d i_N - i_{NP} )\alpha = [i_P, [d,i_N]] \alpha = [[i_N,d],i_P]
\alpha
= [i_N, i_P]_d \alpha \ .
$$
These equalities imply \eqref{trace}.
\end{proof}

The following corollary, a consequence of the compatibility, will be
used in Section \ref{manifold}.
\begin{corollary}\label{corollary}
Let $(P, N)$ be a Poisson-Nijenhuis structure on a manifold. 
For all $f \in
C^{\infty}(M)$,
\begin{equation}\label{identity}
i_{P} (d_Ndf)= - \frac{1}{2} H^P_{I_1}(f),
\end{equation}
where $H^P_{I_1} = Pd  {\rm{Tr}}  N$ is the hamiltonian vector field 
with hamiltonian\break 
$I_1 =  {\rm{Tr}}  N$ in the Poisson structure $P$.
\end{corollary}
\begin{proof}
When $C(P,N) = 0$, formula \eqref{trace} for $\alpha = df $ yields 
\eqref{identity}.
\end{proof}
 
\begin{remark}
When $P$ and $N$ are compatible,
the derived bracket $[i_{N},i_P]_d$ is a
derivation of degree $-1$ of the algebra of forms equal to the
interior product by the vector field $ - \frac{1}{2}
P d  {\rm{Tr}}  N$. A proof of this property can be found in \cite{BM}. 
\end{remark}

\section{The hierarchy of modular classes of a Poisson-Nijenhuis
  manifold}\label{manifold}

\subsection{The modular class of $(TM, N, [ \, , \, ]_N)$.}
Let $N$ be a Nijenhuis tensor on manifold $M$.
Given $\lambda \otimes \mu$, where $\lambda$ is a nowhere vanishing 
multivector of top degree
and $\mu$ a volume element on $M$, the
modular class of the Lie algebroid $(TM, N, [ \, , \, ]_N)$ 
is the class in the
$d_N$-cohomology of the 1-form $\xi^{(N)}$ such that, for all
$X \in \Gamma (TM)$,
$$
< \xi^{(N)} , X > \lambda \otimes \mu = [X , \lambda]_N \otimes \mu +
\lambda \otimes L_{NX} \mu \ .
$$
If $e_1. \ldots. e_n$ is a local basis of $\Gamma(TM)$ such that
$\lambda = e_1 \wedge \ldots \wedge e_n$, then
$$
[X , \lambda]_N = \sum_{j=1}^n (-1)^j [X,e_j]_N e_1 \wedge \ldots \wedge
\widehat{e_j} \wedge \ldots \wedge e_n \ . 
$$
Since $[X,Y]_N = [NX,Y] + (L_XN) Y$, we obtain
$$
[X , \lambda]_N = L_{NX} \lambda + \sum_{j=1}^n (L_{X}N)^j_j e_1 \wedge \ldots \wedge
{e_j} \wedge \ldots \wedge e_n. 
$$
Choosing $\lambda$ and $\mu$ such that $<\lambda , \mu > \, = 1$ which
implies that
$L_{NX} \lambda \otimes \mu + \lambda \otimes L_{NX} \mu =0$, and 
using the relation $ \sum_{j=1}^n (L_{X}N)^j_j =\sum_{j=1}^n
L_{X}(N^j_j)$, we obtain
$$
< \xi^{(N)} , X > \lambda \otimes \mu = i_X (d {\rm{Tr}} N) \, \lambda \otimes \mu \ .
$$
Thus we have recovered the result of \cite{DF}:
\begin{proposition}\label{proposition2.1}
The modular class in the $d_N$-cohomology 
of the Lie algebroid $(TM, N, [ \, , \, ]_N)$ is the
class of 
the $1$-form $d  {\rm{Tr}}  N$.
\end{proposition}

The $d_N$-cocycle $\xi^{(N)}=d {\rm { {\rm{Tr}} }} N$ is in fact independent of the choice of $\lambda
\otimes \mu$. The class it defines 
can also be considered to be the class of the morphism
of Lie algebroids $N \colon 
(TM, N, [ \, , \, ]_N) \to  (TM, {\rm {id}}, [ \, , \, ])$.
 
Similarly, the modular classes associated to the Nijenhuis tensors
$N^k$, $k\in \mathbb N$, $k \geq 2$, are the $d_{N^k}$-classes of the 1-forms 
$d  {\rm{Tr}} (N^k)$.

\subsection{The modular class of a Poisson-Nijenhuis manifold}
We shall now consider the case of a manifold $M$ with a PN-structure.
Let $P_0=P$
and $P_1=NP, \ldots, P_k = N^kP, \ldots$

For each Poisson structure $P_k$ on
$M$, $k \geq 0$, the modular vector field associated to a volume form
$\mu$ 
on $M$ is,
by definition, the $d_{P_k}$-cocycle $X^k_{\mu}$ satisfying
\begin{equation}\label{modular}
< X^k_{\mu} , df > \mu =L_{H^{P_k}_f} \mu \ ,
\end{equation}
for all $f \in C^{\infty}(M)$, that is 
$< X^k_{\mu} , df > \mu = d i_{P_k df} \mu$. 
It follows that, for all 1-forms $\alpha$,
\begin{equation}\label{forms}
< X^k_{\mu} , \alpha > \mu = d i_{P_k \alpha} \mu - (i_{P_k} d\alpha) \mu
\ . 
\end{equation}
 
We now consider the vector fields
\begin{equation}\label{hierarchy}
X^{(k)} = X^{k}_{\mu} - N X^{k-1}_{\mu} \ ,
\end{equation}
for $k \geq 1$, and we list their basic properties:
\begin{itemize}
\item
For each $k$, $X^{(k)}$ is independent of $\mu$. 
It is called the {\it $k$-th modular vector field} of $(M,P,N)$.
\item $X^{(k)}$ is a 
$d_{P_{k}}$-cocycle. Its class is called the {\it $k$-th modular
  class} of the
PN-manifold. In particular, the $d_{NP}$-class of $X^{(1)}$ 
is called {\it the modular class} of $(M,P,N)$.
\item
The $k$-th modular class of $(M,P,N)$ is one-half
the relative modular class of the morphism
of Lie algebroids 
$^tN : (T^*M, P_{k}, [\, , \, ]_{P_{k}})  \to  
(T^*M, P_{k-1}, [\, , \, ]_{P_{k-1}})$.
\end{itemize}

\subsection{Properties of the hierarchy of modular vector fields}
\begin{proposition}
The modular vector fields of a PN-manifold 
$(M,P,N)$ satisfy 
\begin{equation}\label{relation}
 X^{(k)}  = - \frac{1}{2} H^P_{I_{k}} \ , \, k \geq 1, 
\end{equation}
where $I_k = {\rm{Tr}}  \frac{N^k}{k}$, $k \geq 1$, 
is the sequence of fundamental
functions in involution. 
\end{proposition}
\begin{proof}
For clarity, we first prove the case $k=1$.
It follows from formula \eqref{forms} and Corolllary \ref{corollary}
that, for all $f \in C^{\infty}(M)$, 
$$
<NX^0_{\mu}, df > \mu = <X^0_{\mu}, N df > \mu
$$
$$
= di_{PN df} \mu - (i_P d Ndf) \mu 
= di_{NP df} \mu + \frac{1}{2} <P d  {\rm{Tr}}  N, df > \mu \ , 
$$
while 
$$
< X^1_{\mu}, df > \mu
= di_{NP df} \mu \ .
$$
Therefore
$ X^{(1)} = X^1_{\mu} - N X^0_{\mu} = - \frac{1}{2} P d  {\rm{Tr}}  N = -
\frac{1}{2} H^P_{I_1} $.

The case $k \geq 2$ is similar. Applying Corollary \ref{corollary}
to the compatible pair $(N^{k-1}P, N)$, we obtain
$$
<X^{(k)}, df >  = i_{N^{k-1}P} d N df  = i_{N^{k-1}P} d_N df  
= - \frac{1}{2}  <N^{k-1}P d  {\rm{Tr}}  N, df > \ .
$$
The result follows from $N^{k-1}P d  {\rm{Tr}}  N = P N^{k-1} d  {\rm{Tr}}  N
 = P d
{\rm{Tr}}  \frac{N^k}{k}$.
\end{proof}

\begin{remark}
The sequence of modular vector fields $X^{(k)}$, $ k \geq 1$, 
coincides with the well-known
sequence \cite{M} of bihamiltonian vector fields of a PN-manifold.
It follows that 
$X^{(k)} = N X^{(k-1)}$. 
\end{remark}

\begin{remark}
The sequence of modular vector fields of a Poisson-Nijenhuis manifold
introduced by 
Damianou and Fernandes in \cite{DF} is $X_k$, $k \geq 1$, defined by the
recursion $X_1= X_N = X^1_{\mu} - NX^0_{\mu}$ and
$X_k = N X_{k-1}$, for $k \geq 2$. They proved that 
$X_k = - \frac{1}{2} P d {\rm {Tr}} \frac{N^k}{k}$, for $k \geq 1$. 
Though the defintion of
the hierarchy $X^{(k)}$ that we have considered differs from theirs,
the two hierarchies still coincide.

If we denote the modular vector field of the PN-structure $(N,P)$
by $X_{N,P}$, 
then $X^{(k)} = X_{N,N^{k-1}P}$, while $X_k = N^{k-1}X_{N,P}$. 
The vector fields $X_{N,P}$ satisfy
$$
X_{N,NP} +N  X_{N,P} = X_{N^2,P} \ ,
$$
and, more generally, 
$$
X_{N,N^kP} +N  X_{N,N^{k-1}P} = X_{N^2,N^{k-1} P} \ .
$$
This relation is immediate from the definition. Each term is a
hamiltonian vector field with respect to $N^kP$; 
each of the terms on the left-hand side is equal to 
$- \frac{1}{2} P  N^{k} d {\rm {Tr}} N$,
while the right-hand side is 
$- \frac{1}{2} P N^{k-1} d {\rm {Tr}} {N^2} = - P N^{k} d {\rm {Tr}} N$.
\end{remark}

\begin{remark}\label{remark2.3}
It follows from the morphism properties of $P$, $NP$ and $^t N$ that  
the relative modular classes of $P \colon (T^*M, P, [ \, , \,]_P) \to
(TM, Id, [ \, , \,])$,
$NP \colon (T^*M, NP, [ \, , \,]_{NP}) \to
(TM, Id, [ \, , \,])$, and $^t N : 
(T^*M, NP, [ \, , \,]_{NP}) \to (T^*M, P, [ \, , \,]_P)$ are defined and
satisfy
\begin{equation}\label{morphism}
Mod^{NP} - N Mod^P = Mod\,^{^t N} \ .
\end{equation}
A representative of this $d_{NP}$-cohomology class is 
$ - Pd{\rm {Tr}} N
= 2 X^{(1)}$.

More generally, a representative of the modular class 
of the morphism $^t N^{ k }$
from $(T^*M, P_k, [ \, , \, ]_{P_k})$ to $(T^*M, P, [ \, , \, ]_{P})$ is 
 $ - Pd{\rm {Tr}} N^k  =  2 k X^{(k)}$.

\end{remark}

\begin{remark} The modular classes of the morphisms
$N : 
(TM, N, [ \, , \,]_{N}) \to (TM, Id, [ \, , \,])$ 
and $^t N : 
(T^*M, NP, [ \, , \,]_{NP}) \to (T^*M, P, [ \, , \,]_P)$ are related by
\begin{equation}\label{transpose}
Mod \,  ^{^t N} = - P Mod ^N \ .
\end{equation}
\end{remark}
Relation \eqref{transpose} can be generalized in two ways.
\begin{proposition} 
(i)
The modular classes of the morphisms
$$N^{ k } \colon (TM, N^k, [ \, , \, ]_{N^k}) \to (TM, Id, [ \, ,
\, ]) \quad \quad \quad  {\rm{and}} 
$$
$$
^t N^{ k } \colon (T^*M, P_k, [ \, , \, ]_{P_k}) \to (T^*M, P, 
[ \, , \, ]_{P})$$
are related by
$$
Mod \, ^{^t N^{ k }} =- P Mod ^{N^{ k }} \ .
$$
(ii) The modular classes of the morphisms
$$
N^{[ k ]} \colon (TM, N^k, [ \, , \, ]_{N^k}) \to (TM, N^{k-1}, [ \,
, \, ]_{N^{k-1}}) \quad  \quad  \quad  {\rm{and}}
$$ 
$$^t N^{[k]} \colon (T^*M, P_k, [ \, , \, ]_{P_k}) 
\to (T^*M, P_{k-1}, [ \, , \, ]_{P_{k-1}})$$ 
are related by
$$
Mod \,  ^{^t N^{[ k ]}} =- P Mod ^{N^{[ k ]}} \ ,
$$
and a representative of the modular class 
of the morphism $^t N^{[ k ] }$ is $2 X^{(k)}$.
\end{proposition} 
\begin{proof} (i) follows from Proposition \ref{proposition2.1}
and Remark \ref{remark2.3}. To prove (ii), we compute 
a representative of the modular class of $N^{[k]}$,
$$
d {\rm {Tr}} N^k
- \, ^tN d{\rm {Tr}} N^{k -1} =  d{\rm {Tr}} \frac{N^{k}}{k} \ ,$$
and a  
representative of the modular class of $^t N^{[k]}$,
$$ 
2 (X^k_{\mu} - N X^{k-1}_{\mu}) = 2 X^{(k)} =  - Pd{\rm
  {Tr}} \frac{N^k}{k} \ . 
$$
\end{proof}
\begin{remark}
Computations of a representative of $Mod \, ^{^tN^k}$ either from the equality
$2(X^k_{\mu} - N^kX^0_{\mu})= 2
\sum_{\ell = 1}^{k} 
N^{k- \ell} X^{(\ell)}$ or from the equality
$Mod \, ^{^tN^k} \break
= \sum_{\ell = 1}^{k} N^{k-\ell} Mod \,
^{^tN^{[\ell]}}$ both recover the fact, stated in Remark \ref{remark2.3}, that 
a representative of $Mod \, ^{^tN^k}$ is equal to
$ - \sum_{\ell = 1}^{k} N^{k-\ell}  Pd{\rm {Tr}}
\frac{N^{\ell}}{\ell} =  - Pd{\rm
  {Tr}} {N^k} \break
=  2k X^{(k)}$. 
\end{remark}

\bigskip

\noindent {\it Acknowledgments} We thank
Rui Fernandes and Camille Laurent-Gengoux 
for useful exchanges on the topic
of this letter.

\bigskip

\hfill Centre de Math\'ematiques Laurent Schwartz

\hfill Ecole Polytechnique, F-91128 Palaiseau, France

\hfill yks@math.polytechnique.fr

\bigskip

\hfill Department of Mathematics

\hfill University of Milano Bicocca

\hfill Via Corsi 58, I-20126 Milano, Italy

\hfill magri@matapp.unimib.it

\end{document}